\documentclass[12pt, reqno]{amsart}
\usepackage{amsfonts}
\usepackage{bbm}
\usepackage{amscd,amsfonts}
\usepackage{amssymb, eucal, amsfonts, amsmath, xypic, latexsym, tikz}
\usepackage{pifont}
\usepackage{mathrsfs,color}
\usepackage{amsthm,indentfirst,bm,fancyhdr,dsfont}
\usepackage{graphicx}
\usepackage[all]{xy}
\usepackage[CJKbookmarks=true]{hyperref}
\usepackage[normalem]{ulem}
\usepackage{youngtab}

\usepackage{mathrsfs}
\usepackage{amsmath}
\usepackage{amssymb}
\usepackage{hyperref}

\newtheorem{theorem}{Theorem}[section]
\newtheorem{lemma}[theorem]{Lemma}

\newtheorem{proposition}[theorem]{Proposition}
\newtheorem{corollary}[theorem]{Corollary}
\theoremstyle{definition}

\newtheorem{defn}[theorem]{Definition}
\newtheorem{definition}[theorem]{Definition}

\newtheorem{remark}[theorem]{Remark}

\numberwithin{equation}{section}

\newcommand{\pa}[1]{|{#1}|}
\def\ggg{\mathfrak{g}}

\def\gl{\mathfrak{gl}}

\def\g{\mathfrak{g}}
\def\ggg{\mathfrak{g}}

\def\mmm{\mathfrak{m}}
\def\ppp{\mathfrak{p}}

\def\ttt{\mathfrak{t}}
\def\hhh{\mathfrak{h}}

\def\bbc{\mathbb{C}}

\def\bz{{\bar 0}}

\def\ad{\mathsf{ad}}

\newcommand{\tp}{\operatorname{pa}}

{\vskip-\lastskip\medskip
 \noindent
 }%
{\qed\par\medskip
 }

\def\sfs{\textsf{s}}

\begin{document}

\title[On the one-dimensional representations of finite $W$-superalgebras for $\gl_{M|N}$]
{On the one-dimensional representations of finite $W$-superalgebras for $\gl_{M|N}$
}

\author{Fanlei Yang and Yang Zeng}

\begin{center}
\textsc{Dedicated to Professor Bin Shu on the occasion of his 60th birthday}
\end{center}

\address{School of Mathematical Sciences,  East China Normal University, No. 500 Dongchuan Rd., Shanghai 200241, China} \email{52275500005@stu.ecnu.edu.cn}

\address{School of Mathematics, Nanjing Audit University, No. 86 West Yushan Rd., Nanjing, Jiangsu Province 211815, China}
\email{zengyang@nau.edu.cn}

  \thanks{This work is partially supported by the National Natural Science Foundation of China (Nos. 12461005, 11701284).}

\begin{abstract} Let $\ggg=\gl_{M|N}(\Bbbk)$ be the general linear Lie superalgebra over an algebraically closed field $\Bbbk$ of characteristic zero. Fix an arbitrary even nilpotent element $e$ in $\ggg$ and let $U(\ggg,e)$ be the finite $W$-superalgebra associated to the pair $(\ggg, e)$. In this paper we will  give a complete classification of one-dimensional representations for $U(\ggg,e)$. To achieve this, we use the tool of shifted super Yangians to determine the commutative quotients of the finite $W$-superalgebras.
\end{abstract}
\maketitle

\section{Introduction}


\subsection{}\label{1.1}
A finite $W$-algebra $U(\mathfrak{l},e)$ is a certain associative algebra associated with a complex semi-simple Lie algebra  $\mathfrak{l}$ and a nilpotent element $e\in\mathfrak{l}$. The study of finite $W$-algebras can be traced back to Kostant's work in the case when $e$ is regular \cite{Ko}. Along the proof of the celebrated Kac-Weisfeiler conjecture for Lie algebras of reductive groups in \cite{P1}, Premet further developed the finite $W$-algebras in full generality in \cite{P2} from the viewpoint of modular representation theory. Furthermore, finite $W$-algebras can be interpreted as the quantizations of Slodowy slices \cite{GG, P2}. Since then, finite $W$-algebras have appeared in many branches of mathematics. We refer the readers to some survey papers (for example, \cite{A, L4, W}) for details.


For the special case with $\mathfrak{l}=\gl_N$, the connection between Yangians and finite $W$-algebras associated to an arbitrary nilpotent $e\in\gl_N$ was established by Brundan-Kleshchev in \cite{BK06}. The main result of \cite[Theorem 10.1]{BK06} can be summarized as follows: the finite $W$-algebras over the field of complex numbers associated with a nilpotent $e\in\gl_N$ are isomorphic to a quotient of the shifted Yangians associated to $\gl_n$, where $n$ is the number of Jordan blocks of $e$. Moreover, an explicit realization of finite $W$-algebra of type $A$ in terms of the generators and their relations is obtained. This provides a powerful tool for the study of their representation theory as discussed in \cite{BGK, BK08, GT}.
\subsection{}
Among all the representations of finite $W$-algebras, the one-dimensional representations play an especially important role. Let us list some related consequences as below.  Firstly, the images under the Skryabin's equivalence as in \cite{Sk} are all completely prime, and therefore play a key role in Joseph's theory of Goldie rank polynomials (see \cite{Lo2}). Secondly, for a complex connected reductive algebraic group $L$ with Lie algebra $\mathfrak{l}$, those one-dimensional representations classify quantizations of $L$-equivariant coverings of nilpotent $L$-orbits $\mathcal{O}$ (see \cite{Lo1}). Thirdly, Topley studied the variety of one-dimensional representations of a finite $W$-algebra, giving a precise description of the dimensions of the irreducible components. And then this result is applied to prove a conjecture of Losev describing the image of his orbit method map (see \cite{T}).
Finally, the result on the existence of one-dimensional representations of finite $W$-algebras can directly 
lead to the accessibility of lower-bounds of dimensions in the modular representations of reductive Lie algebras predicted by the afore-mentioned Kac-Weisfeiler conjecture, which is known as Humphreys' conjecture (see \cite{Pre, PT2}).


The existence of one-dimensional representations for all finite $W$-algebras over $\mathbb{C}$ was conjectured by Premet in \cite{P3}.  This conjecture was reduced to the rigid case in \cite{Pre}, and then confirmed for classical cases in \cite{PT1}. 
The case of exceptional Lie algebras was dealt with in \cite{GRU,P9}.  The arguments in \cite{P9,PT1} relied on the results proved by Losev in \cite{Lo1,Lo} whilst \cite{P9} was entirely based on computer-aided computations.

Now we turn back to the special case with $\mathfrak{l}=\gl_N$. As mentioned in \S\ref{1.1}, the generators and their relations of corresponding finite $W$-algebras are introduced by Brundan-Kleshchev in \cite{BK06}.
In virtue of these explicit presentations, Premet described the maximal commutative quotient of finite $W$-algebras in \cite[Theorem 3.3]{Pre}, which helps Brundan to obtain a classification of one-dimensional representations of these finite $W$-algebras in \cite[\S2]{Bru}. Following similar methods of above, the modular counterpart was achieved by Goodwin-Topley in \cite[\S5]{GT}.

\subsection{}\label{wsuperin} In the same time, the theory of finite $W$-superalgebras was also developed. Given a basic classical Lie superalgebra $\ggg$, the finite $W$-superalgebras can be defined in a very similar way as the Lie algebra case, except that the nilpotent element $e\in\ggg$ is assumed to be even with other modifications.

In the work of Wang-Zhao \cite{WZ}, they initiated the study of modular representations of basic classical Lie superalgebras in positive characteristic, formulating the super Kac-Weisfeiler property for those Lie superalgebras as well as presenting the definition of modular $W$-superalgebras. As for the finite $W$-superalgebras over $\bbc$,  Zeng-Shu gave the PBW structure theorem in \cite{ZS1}, which shows that the construction of finite $W$-superalgebras can be divided into two cases (the even case and the odd case) in virtue of the parity of the dimension for a certain subspace of $\ggg$,  thus makes the construction of finite $W$-superalgebras significantly different from that of finite $W$-algebras.
As a further work, Zeng-Shu studied the modular representations of basic Lie superalgebras in \cite{ZS4}, as a remarkable application of finite $W$-superalgebras. Apart from the above, the topics on finite $W$-superalgebras and their representations have
been extensively studied from various aspects in recent years, and we refer the readers to \cite{BBG2,CC,P,PS2,PS3,SX,suh,W,Z2} etc. for more details.

For the connections between finite $W$-superalgebras and super Yangians, Briot-Ragoucy obtained the very first result in \cite{BR}, saying that if the nilpotent element $e\in\gl_{M|N}$ is rectangular, then the associated finite $W$-superalgebra is isomorphic to a certain quotient of the truncated super Yangian.
From then on, there have been some results \cite{BBG2, Pe2, Pe3} generalizing the above observation when the nilpotent element $e$ satisfies some assumptions. Finally, Peng \cite{P} has established a superalgebra isomorphism between the finite $W$-superalgebras associated with an arbitrary even nilpotent element $e\in\gl_{M|N}$ and a quotient of the shifted super Yangians, obtaining a super analogue of the main result of \cite{BK06} for Lie superalgebras of type $A$ in full generality. 
The connection between the finite $W$-superalgebra for the queer Lie superalgebras and the super Yangians was considered by Poletaeva-Serganova in \cite{PS2,PS3}.

\subsection{} Parallel to the finite $W$-algebra case, it is a natural to consider the representation theory of finite $W$-superalgebras. For the finite $W$-superalgebras attached to the principal nilpotent orbits in $\gl_{M|N}$, their finite-dimensional modules are considered by Brundan-Goodwin in \cite{BG3}. 

For the topics related to one-dimensional representations of finite $W$-superalgebras, there are already some results on the existence of these modules. It is conformed by Zeng-Shu in the cases for $\gl_{M|N}$ with artbitrary nilpotent orbits in \cite[Proposition 4.7]{ZS4}. For the situation associated with $\mathfrak{osp}_{1|2n}$ and  principal nilpotent orbits, the corresponding finite $W$-superalgebras belong to the odd case as aforementioned in \S\ref{wsuperin}. By virtue of \cite[Proposition 3.7]{ZS4}, those finite $W$-superalgebras do not have one-dimensional representations. As an alternative consequence, they admits two-dimensional representations by \cite[Proposition 5.8]{ZS4}.
As a more general result,
Zeng-Shu \cite{ZS5} showed that the finite $W$-superalgebras attached to any basic classical Lie superalgebras and minimal  nilpotent orbits admit one-dimensional or two-dimensional representations, depending on which case finite $W$-superalgebras belongs to. 
Apart from all above, the classification of those one-dimensional modules and their relation to other aspects of the representation theory of Lie superalgebra remain to be open.

\subsection{}
In the paper we will concentrate on the general linear Lie superalgebra attached to arbitrary nilpotent orbit.
Let $\ggg=\gl_{M|N}$, and $\ppp$ be an appropriate parabolic subalgebra of $\ggg$ given in \eqref{pppmmm}. The goal of this paper is to generally establish the classification of one-dimensional representations for finite $W$-superalgebras for type $A$. 
That is, we explicitly give all possible one-dimensional modules for $U(\ppp)$, which can be induced to obtain the one-dimensional representations for finite $W$-superalgebras, as a super analogue of the result of \cite{Bru,GT} for Lie algebras of type $A$ in full generality.

We shortly explain our approach, which is basically generalizing the arguments in \cite{GT}  to the general linear Lie superalgebras with suitable modifications. Although there are similarities between $\gl_N$ and $\gl_{M|N}$ and also their associated finite $W$-(super)algebras, some of the earlier approaches are no longer available in our case.

As the first step, we introduce the definition of a pyramid $\pi$ and the associated finite $W$-superalgebra $U(\ggg,e)$. Then the definitions of  column-connected associated with a $\pi$-tableau $\text{Tab}_\Bbbk(\pi)$ (which is strikingly different from the original one in the non-super case) and row-equivalent between two tableaux are given. For any $A\in\text{Tab}_\Bbbk(\pi)$, the weights of $A$ are defined, which can help us to construct one-dimensional modules for $U(\ppp)$. Then by inducing this module we can obtain one-dimensional representations of $U(\ggg,e)$, which is the subalgebra of $U(\ppp)$.

We further introduce the shifted super Yangian $Y_{m|n}(\sigma)$ of $\gl_{M|N}$ associated to a pyramid $\pi$, and then define the truncated shifted super Yangian $Y_{m|n}^\ell(\sigma)$ of level $\ell$. Let $U({\mathfrak{g}},e)^{\text{ab}}$ denote the commutative quotient of the finite $W$-superalgebra $U({\mathfrak{g}},e)$ as defined in \S\ref{W_pi}.
With the help of the commutative quotient of $Y_{m|n}^\ell(\sigma)$, we obtain the first main result of the paper as below.
\begin{theorem}(see Theorem \ref{glmn})
The commutative quotient $U({\mathfrak{g}},e)^{\text{ab}}$ is isomorphic to a polynomial algebra in $\ell=p_{m+n}$ variables \begin{align*}
\big\lbrace d_{i}^{(r)} \,|\, {i=1,\ldots,m+n,\;  1\leq r\leq p_i-p_{i-1}}\big\rbrace,
\end{align*}where the meaning of the $d_{i}^{(r)}$'s is given in Proposition \ref{ab spanning set}.
\end{theorem}

For the next step, we obtained a consequence on the elementary symmetric polynomials, and then state the second main result of the paper, for which the consequence in Theorem \ref{glmn} is applied in its proof.
\begin{theorem}(see Theorem \ref{mean})
For any $A\in\text{Tab}_\Bbbk(\pi)$, there is a one-dimensional $U(\ggg,e)$-module $\widetilde{\Bbbk}_{\overline{A}}$, which extends the action of $U(\ggg,e)^0$ on $\Bbbk_{\overline{A}}$ if and only if $A$ is row-equivalent to a column-connected tableau, where $U(\ggg,e)^0$ is a proper subalgebra $U(\ggg,e)$ and $\Bbbk_{\overline{A}}$ is a one-dimensional $U(\ggg,e)^0$-module defined as in \S\ref{one-dimensional modules for U(g,e)}.
\end{theorem}

In the final part, we will give an application of the above theorem. 
For the 
even non-degenerate supersymmetric invariant bilinear form $(\cdot,\cdot)$ on $\ggg=\ggg_{\bar 0}\oplus\ggg_{\bar 1}=\gl_{M|N}(\Bbbk)$ given in \eqref{strace},
define $\chi\in{\ggg}^{*}$ by letting $\chi(x)=(e,x)$ for all $x\in{\ggg}$.
Denote by $\ggg^e$ for the centralizer of $e$ in $\ggg$, and write $d_0=\text{dim}\,\ggg_{\bar 0}-\text{dim}\,\ggg^e_{\bar 0}$, $d_1:=\text{dim}\,\ggg_{\bar 1}-\text{dim}\,\ggg^e_{\bar 1}$. 

Let $\ggg_{\mathbb{F}}$ be the modular counterpart of $\ggg$, and still write $\chi\in({\ggg}_{\mathbb{F}})^*_{\bar0}$ for the modular version of $\chi\in{\ggg}^{*}$ as defined above. Then $\chi$ is a nilpotent $p$-character. 
Denote by $U_\chi({\ggg}_{\mathbb{F}})$ the reduced enveloping algebra of $\ggg_{\mathbb{F}}$ associated to $\chi$. 

As a direct corollary of Theorem \ref{mean}, we will provide another proof of the super version of Humphreys' conjecture for $\gl_{M|N}(\mathbb{F})$ in characteristic $p\gg0$, for which the original one is obtained by Zeng-Shu in \cite{ZS2}.
\begin{corollary}\label{mincor}
Let $\ggg_\mathbb{F}=\mathfrak{gl}_{M|N}(\mathbb{F})$ over a field $\mathbb{F}$ of characteristic $p\gg0$. For any nilpotent $p$-character $\chi\in\ggg^*_{\bar0}$, the reduced enveloping algebra $U_\chi(\ggg_\mathbb{F})$ admits irreducible modules of dimension $p^{\frac{d_0}{2}}2^{\frac{d_1}{2}}$.
\end{corollary}


\subsection{} The paper is organized as follows. In \textsection \ref{pyramids}, we first recall the definition of finite $W$-superalgebra $U(\ggg,e)$ and its associated pyramid, and also the tableaux and weights associated to them. Then we classify the one-dimensional modules for an appropriate parabolic subalgebra $\ppp$ (see \eqref{pppmmm}) of $\ggg$, which helps us to establish the main result in \textsection \ref{One-dimensional modules}. In \textsection \ref{Truncated shifted super Yangian}, we recall some necessary facts about the truncated shifted super Yangian $Y_{m|n}^\ell(\sigma)$ and also the isomorphism between $U(\ggg,e)$ and $Y_{m|n}^\ell(\sigma)$. In virtue of these, we can determine the structure of their commutative quotients. 
In the concluding section \textsection \ref{One-dimensional modules}, we finally obtain the classification of one-dimensional representations for finite $W$-superalgebras. And then we will give another proof of \cite[Proposition 3.1]{ZS2}, which assures the existence of the representations of $\ggg_\mathbb{F}$ of minimal dimension given in Corollary \ref{mincor}.
\subsection{}
Throughout we work with an algebraically closed field $\Bbbk$ of characteristic zero, or the algebraically closed field $\mathbb{F}=\overline{\mathbb{F}}_p$ of positive characteristic $p\gg0$ as the ground field.
We consider vector spaces, subalgebras, ideals, modules, and submodules, {\textit{etc}}. in the super sense unless otherwise specified, throughout the paper. 

The isomorphism between superalgebras is assumed to be a $\mathbb{Z}_2$-graded parity-preserving linear map that is an isomorphism in the usual sense.

\section{Finite $W$-superalgebras and pyramids}\label{pyramids}
In this section, we first recall the definition of a finite $W$-superalgebra $U(\gl_{M|N},e)$, which is determined by an even nilpotent element $e$ of $\ggg=\gl_{M|N}$. And then we give the one-dimensional modules for an appropriate parabolic subalgebra $\ppp$ of $\ggg$, which plays a key role in the understanding of the corresponding $U(\ggg,e)$-modules. To achieve these, we will use the tool of pyramids and tableaux.

Based on \cite{P}, the contents of this section is roughly a generalization of the finite $W$-algebra case as in \cite[\S2]{Bru} and \cite[\S2]{GT}, with a lot of modifications. One can observe that the emergence of odd parts in $\ggg$ makes the situation much more complicated.

\subsection{Finite $W$-superalgebras of $\gl_{M|N}$}\label{Finite $W$-superalgebras}
In this subsection we will recall some knowledge of the finite $W$-superalgebras.
Let $\mathfrak{g}=\gl_{M|N}$ be the general linear Lie superalgebra over $\Bbbk$, and
we can identify $\mathfrak{g}$ with the set of $(M+N)\times(M+N)$ matrices with the standard $\mathbb{Z}_2$-grading $\ggg=\ggg_{\bar 0}\oplus \g_{\bar 1}$. Let $\mathfrak{t}$ be a standard Cartan subalgebra of ${\ggg}$ consisting of diagonal matrices, and let $\Phi$ be the root system of ${\ggg}$ relative to $\mathfrak{t}$.

Let $(\cdot,\cdot)$ be the non-degenerate supersymmetric even invariant bilinear form on $\ggg$ defined by
\begin{equation}\label{strace}
(x,y):=\text{str}(xy)
\end{equation}
for all $x,y \in \ggg$, where 
str denotes the supertrace. Every elements of $\ggg$ appearing in any equations are assumed to be $\mathbb{Z}_2$-homogeneous unless specifically mentioned.
\begin{definition}Let $e$ be a nilpotent element in $\ggg_{\bar{0}}$ and $h$ be a semisimple element.
We say $(e,h)$ is a {\it good pair} if $\ad\,h$ gives an even good $\mathbb{Z}$-grading of $\ggg$ for $e$, which means the following conditions are satisfied:
\begin{itemize}
\item[(1)] $\ad\,h(e)=2e$,
\item[(2)] $\ggg=\bigoplus_{j\in\mathbb{Z}}\ggg(j)$, where $\ggg(j):=\{x\in \ggg \mid \ad h(x)=jx\}$,
\item[(3)] the center of $\ggg$ is contained in $\ggg(0)$,
\item[(4)] $\ad\,e:\ggg(j)\rightarrow \ggg(j+2)$ is injective for all $j\leq -1$,
\item[(5)] $\ad\,e:\ggg(j)\rightarrow \ggg(j+2)$ is surjective for all $j\geq -1$.
\item[(6)] the $\mathbb{Z}$-grading is even; that is, $\ggg(i)=0$ for all $i\notin 2\mathbb{Z}$.
\end{itemize}
\end{definition}

In fact, for any even nilpotent $e\in\gl_{M|N}$, by \cite[Theorem 2.4]{P} we can always find some $h$ 
such that $(e,h)$ is a good pair (see \S\ref{Pyramids and W-superalgebras} for more details).
We can further assume that $h\in\mathfrak{t}$.

Fix a good pair $(e,h)$ in $\ggg$. Define the following subalgebras of $\ggg$ by
\begin{equation}\label{pppmmm}
\ppp:=\bigoplus_{j\geq 0}\ggg(j),\qquad\hhh:=\ggg(0),\qquad\text{and}\quad\mmm:=\bigoplus_{j<0}\ggg(j).
\end{equation}
Define $\chi \in \ggg^*$ by
$$\chi(y):=(e,y),\qquad\forall y\in \ggg.$$
The restriction of $\chi$ on $\mmm$ extends to a one-dimensional $U(\mmm)$-module. Let $I_\chi$ be the left ideal of $U(\ggg)$ generated by
$$\{ a-\chi(a)\mid a\in \mmm\}.$$
We have $U(\ggg)=I_\chi\oplus U(\ppp)$ and then the following identification
$$U(\ggg)/I_\chi \cong U(\ppp),$$
by the natural projection $\text{pr}_\chi:U(\ggg)\rightarrow U(\ppp).$ One can define the following $\chi$-twisted action of $\mmm$ on $U(\ppp)$ by
$$a\cdot y:=\text{pr}_\chi([a,y])$$
for all $a\in\mmm,y\in U(\ppp).$

\begin{defn}\label{eq: The finite $W$-superalgebra}
The finite $W$-superalgebra $U(\ggg,e)$ is defined to be the space of $\mmm$-invariants in $U(\ppp)$ under the $\chi$-twisted action; or to say,
\begin{align*}
U(\ggg,e):=U(\ppp)^{\mmm}&=\{y\in U(\ppp)\mid \text{pr}_\chi([a,y])=0,\forall a\in\mmm\} \cr
&=\{y\in U(\ppp)\mid (a-\chi(a))y\in I_\chi,\forall a\in\mmm\}.
\end{align*}
\end{defn}

\subsection{Pyramids and finite $W$-superalgebras}\label{Pyramids and W-superalgebras}
The contents of this subsection is based on \cite{P}, and more details can be seen there.
\subsubsection{}\label{basic notations}
We fix a partition $\mathbf{p}=(p_1,\ldots,p_{m+n})$ on $M+N$ with $p_1\leq \ldots\leq p_{m+n}$. A pyramid $\pi$ associated with $\mathbf{p}$ is a diagram with $p_{m+n}$ boxes in the bottom row, $p_{m+n-1}$ boxes in the row above it, and so on, stacked in such way that any box not in the bottom row, lies directly above a box below it in the row
beneath it, and the boxes occur consecutively in each row.


Put ``$+$" or ``$-$" in the boxes of a pyramid such that every boxes in a row have the same ``$+$" or ``$-$" labeling. Assume that we have
$M$ (resp. $N$) boxes labeled with ``$+$" (resp. ``$-$"), and also $m$ (resp. $n$) rows labeled by ``$+$" (resp. ``$-$") in $\pi$.
We enumerate those ``$+$" boxes by $1,2,\ldots,M$ down columns from left to right, and enumerate those ``$-$" boxes by $\overline{1},\overline{2},\ldots,\overline{N}$ also in this way.

Write $\ell=p_{m+n}$.
The columns of $\pi$ are labelled $1,2,\ldots,\ell$ from left to right, and the rows are labelled $1,2,\ldots,m+n$ from top to bottom. We write $\text{row}(i)$ and $\text{col}(i)$ for the row and column of the $i$-th box, respectively.

A $0^m1^n$-sequence, or $01$-sequence for short, is an ordered sequence $\Upsilon$ consisting of $m$ $0$'s and $n$ $1$'s. Associated to $\pi$, we can obtain a $0^m1^n$-sequence $\Upsilon$ by assigning the $i$-th digit of $\Upsilon$ to be 0 (resp. 1) if the boxes in the $i$-th row are labeled by $``+"$ (resp. $``-"$).
For $1\leq i\leq m+n$, the $i$-th digit of $\Upsilon$ is denoted by $|i|$.

We assume that each box of $\pi$ is of size $2\times2$ and our pyramid is built on the $x$-axis, with the center of $\pi$ being exactly located above the origin. Let $\text{col}_x(i)$ denote the $x$-coordinate of the center of the box numbered with $i\in I$.
For example, when $\ggg=\gl_{3|6}$, we can construct a pyramid as follows:
\begin{equation}\label{exp}
\pi={\begin{picture}(-70, 70)%
\put(-220,-10){\line(1,0){80}}
\put(-220,10){\line(1,0){80}}
\put(-200,30){\line(1,0){60}}
\put(-200,50){\line(1,0){40}}
\put(-220,-10){\line(0,1){20}}
\put(-200,-10){\line(0,1){60}}
\put(-180,-10){\line(0,1){60}}
\put(-160,-10){\line(0,1){60}}
\put(-140,-10){\line(0,1){40}}
\put(-215,-3){$-$}\put(-195,-3){$-$}\put(-175,-3){$-$}\put(-155,-3){$-$}
\put(-195,16){$+$}\put(-175,16){$+$}\put(-155,16){$+$}
\put(-195,35){$-$}\put(-175,35){$-$}
\put(-80,-35){$x$-coordinates:}
\put(15,-10){\line(1,0){80}}
\put(15,10){\line(1,0){80}}
\put(35,30){\line(1,0){60}}
\put(35,50){\line(1,0){40}}
\put(15,-10){\line(0,1){20}}
\put(35,-10){\line(0,1){60}}
\put(55,-10){\line(0,1){60}}
\put(75,-10){\line(0,1){60}}
\put(95,-10){\line(0,1){40}}
\put(23,-3){$\overline{1}$}\put(43,-3){$\overline{3}$}\put(63,-3){$\overline{5}$}\put(83,-3){$\overline{6}$}
\put(43,16){$1$}\put(63,16){$2$}\put(83,16){$3$}
\put(43,35){$\overline{2}$}\put(63,35){$\overline{4}$}
\put(0,-20){\line(1,0){110}}
\put(53,-23){$\bullet$}
\put(63,-35){$1$}\put(83,-35){$3$}
\put(35,-35){$-1$}\put(15,-35){$-3$}
\end{picture}}\\[15mm]
\end{equation}
We have that \begin{center}
$m=1, n=2, \mathbf{p}=(2,3,4), \Upsilon=101$,\\
$\text{col}_x(1)=-1, \text{col}_x(2)=1, \text{col}_x(3)=3$,\\
$\text{col}_x(\overline1)=-3, \text{col}_x(\overline2)=-1, \text{col}_x(\overline3)=-1, \text{col}_x(\overline4)=1, \text{col}_x(\overline5)=1, \text{col}_x(\overline6)=3.$
\end{center}

Set $h=m-n$ and let $(\check{q}_1,\ldots,\check{q}_\ell)$ denote the {\em super column heights} of $\pi$, where each $\check{q_i}$ is defined to be the number of boxes labeled with ``$+$" subtract the number of boxes labeled with ``$-$" in the $i$-column of $\pi$. Let $\check{\text{row}}(i)$ denote the {\em super row number of the $i$-th box} which is defined to be the number of rows labeled with ``$+$" subtract the number of rows labeled with ``$-$", counting 
from the top row to the row where the $i$-th box lies in.
For example, for the pyramid \eqref{exp} we can obtain
\begin{center}
$\check{q}_1=-1,\check{q}_2=-1,\check{q}_3=-1,\check{q}_4=0,$\\
$\check{\text{row}}(1)=\check{\text{row}}(2)=\check{\text{row}}(3)=0,$\\
$\check{\text{row}}(\overline{1})=\check{\text{row}}(\overline{2})=\check{\text{row}}(\overline{3})=\check{\text{row}}(\overline{4})=\check{\text{row}}(\overline{5})=\check{\text{row}}(\overline{6})=-1.$
\end{center}

From $\pi$, we can define the shift matrix $\sigma=(s_{i,j})$ as follows. For $1\leq i<j\leq m+n$, we let $s_{j,i}$ be the left indentation of the $i$-th row of $\pi$ relative to the $j$-th row,  and we let $s_{i,j}$ be the right indentation of the $i$-th row of $\pi$ relative to the $j$-th row, and set $s_{i,i}=0$. Then we have $p_i=\ell-s_{i,m+n}-s_{m+n,i}$. For example, associated with $\pi$ in \eqref{exp} we have the following shift matrix
\begin{equation*}\label{sigex}
\sigma = \left(\begin{array}{lll}
0&1&1\\
0&0&0\\
1&1&0
\end{array}\right).
\end{equation*}

\subsubsection{}Now we can obtain a good pair $(e,h)$ from a certain pyramid $\pi$ with $M+N$ boxes.
Fix an ordered index set $I=\{1<\ldots<M<\overline{1}<\ldots<\overline{N}\}$, and let $\{v_i\mid i\in I\}$ be the standard basis of $\Bbbk^{M|N}$ with respect to the order
\[v_i<v_j \text{ if } i<j \text{ in } I.\]
Let $\tp(i):=0$ if $i\in\{1,\ldots,M\}$ and $\tp(i):=1$ otherwise.
Let $\{e_{i,j}\,|\, i,j\in I\}$ denote the elementary matrices in $\gl_{M|N}$ with the parity $(\tp(i)+\tp(j)) (\text{mod}\,2)$.
Define the element
\begin{equation}\label{defineepi}
e_\pi:= \sum_{\substack{\text{row}(i) = \text{row}(j) \\ \text{col}(i) = \text{col}(j) - 1}} e_{i,j} \in\ggg_{\overline{0}}. 
\end{equation}
Define the diagonal matrix
\begin{equation*}\label{definehpi}
h_\pi:=-\text{diag}(\text{col}_x(1),\ldots,\text{col}_x(M),\text{col}_x(\overline{1}),\ldots,\text{col}_x(\overline{N})).
\end{equation*}
Then one can easily check that $(e_\pi,h_\pi)$ is a good pair.
For example, the elements $e_\pi$ and $h_\pi$ associated with the pyramid in (\ref{exp}) are
\begin{equation*}
\begin{split}
e_\pi&=e_{1,2}+e_{2,3}+e_{\overline{2},\overline{4}}+e_{\overline{1},\overline{3}}+e_{\overline{3},\overline{5}}+e_{\overline{5},\overline{6}},\\
h_\pi&=\text{diag}(1,-1,-3,3,1,1,-1,-1,-3).
\end{split}
\end{equation*}

Associated with the good pair $(e_\pi,h_\pi)$ as above, we can obtain a good $\mathbb{Z}$-grading of $\ggg$:
\begin{align}\label{good grading}
\ggg=\bigoplus_{k\in\mathbb{Z}}\ggg(k),\quad\text{where }  \ggg(k)&=\{x\in \ggg \mid \ad\,h_\pi(x)=kx\} \cr
&=\text{span}\{e_{i,j}\in\ggg\mid \text{col}_x(j)-\text{col}_x(i)=k\}.
\end{align}
Under the above settings, the subalgebras $\ppp$ and $\mmm$ can be defined as in \eqref{pppmmm}. So we can introduce the corresponding finite $W$-superalgebra $U(\ggg,e)$ in Definition \ref{eq: The finite $W$-superalgebra}.


%

\subsection{Tableaux and weights}\label{Tableaux}
This subsection is devoted to the discussion on the tableaux and weights, which is a generalization of the non-super case. 
One can observe significant difference in the definition of the column-connected tableau between the non-super version as in \cite[\S2.4]{GT} and the one we give in Definition \ref{cc}. Such a distinction will be explained in the forthcoming Proposition \ref{one dim to column-connected}.
\subsubsection{}
A $\pi$-tableau is a diagram obtained by filling the boxes of $\pi$ with elements of $\Bbbk$. Denote by $\text{Tab}_\Bbbk(\pi)$ the set of all tableaux of shape $\pi$. For $A\in \text{Tab}_\Bbbk(\pi)$, we write $a_i$ for the entry in the $i$-th box of $A$.
\begin{definition}\label{rowe}
Two tableaux are called {\it row-equivalent} if one can be obtained from the other by permuting the entries in the rows.
\end{definition}


\begin{definition}\label{cc}
A tableau $A\in \text{Tab}_\Bbbk(\pi)$ is {\it column-connected} if whenever the $j$-th box of $\pi$ is directly below the $i$-th box and $\tp(i)=\tp(j)$ we have $a_i=a_j+1$; whenever the $j$-th box of $\pi$ is directly below the $i$-th box and $\tp(i)\neq\tp(j)$ we have $a_i+a_j=-1.$
\end{definition}

Let
$\{\varepsilon_1,\ldots,\varepsilon_M,\varepsilon_{\overline1},\ldots,\varepsilon_{\overline N}\}$ be the standard basis of $\ttt^*$
defined by $\varepsilon_i(\text{diag}(d_1,\dots,d_n)) = d_i$.
For $A\in \text{Tab}_\mathbb{\Bbbk}(\pi)$, we define a weight $\lambda_A \in \ttt^*$ by
\begin{equation}\label{lambda_A}
\lambda_A:=\sum\limits_{i\in I}a_i\varepsilon_i.
\end{equation}

To understand the required weights, we will give a decomposition of $\Phi$.
For any $\alpha\in\Phi$, we denote its parity by $\tp(\alpha)$.
We define
\begin{align*}
\Phi^+ &:= \{\varepsilon_i - \varepsilon_j \in \Phi \mid \text{row}(i) <\text{row}(j)\}, \\
\Phi^0 &:= \{\varepsilon_i - \varepsilon_j \in \Phi \mid \text{row}(i) = \text{row}(j)\}, \\
\Phi^- &:= \{\varepsilon_i - \varepsilon_j \in \Phi \mid \text{row}(i) > \text{row}(j)\}.
\end{align*}
Also set
\begin{align*}
\Phi(+) &:= \{\varepsilon_i - \varepsilon_j \in \Phi \mid \text{col}(i) < \text{col}(j)\}, \\
\Phi(0) &:= \{\varepsilon_i - \varepsilon_j \in \Phi \mid \text{col}(i) = \text{col}(j)\}, \\
\Phi(-) &:= \{\varepsilon_i - \varepsilon_j \in \Phi \mid \text{col}(i) > \text{col}(j)\}.
\end{align*}
Then for $\eta,\xi \in \{-,0,+\}$, we define
$$
\Phi(\eta)^\xi = \Phi(\eta) \cap \Phi^\xi.
$$
Note that $\Phi(0)^0 = \varnothing$ and $\Phi^+ \cup \Phi(+)^0$ is the system of positive
roots corresponding to a Borel subalgebra $\mathfrak{b}$ of $\ggg$. Furthermore, $\Phi(+) \cup  \Phi(0)^+$ is
the system of positive roots corresponding to a Borel subalgebra
contained in $\ppp$, and $\Phi(-) \cup \Phi(0)^+$ is another choice of a system of positive roots.

\subsubsection{}
We define
\begin{align*}\label{bar{rho}}
\overline{\rho}:=-\sum\limits_{i\in I}(-1)^{\tp(i)}((\check{q}_1+\ldots+\check{q}_{\text{col}(i)-1})+\check{\text{row}}(i)-(h-\check{q}_{\text{col}(i)}))\varepsilon_i,
\end{align*}
which is a shifted half algebraic sum of positive roots for $\Phi(+) \cup  \Phi(0)^+$ such that
\begin{equation*} \label{e:barrho}
\bar \rho =    \frac{1}{2} \left(\sum_{\alpha \in \Phi(+) \cup \Phi(0)^+}(-1)^{\tp(\alpha)}\alpha\right) - \delta,
\end{equation*}
where
$$
\delta = \frac{M-N+1}{2}\sum_{i\in I}\varepsilon_i.
$$

Define
\begin{align}\label{eta}
\eta:=\sum\limits_{i\in I}(-1)^{\tp(i)}(h-\check{q}_{\text{col}(i)}-\cdots-\check{q}_\ell)\varepsilon_i,
\end{align}
and
\begin{equation*}\label{rho_h}
\rho_\hhh:=-\sum\limits_{i\in I}(-1)^{\tp(i)}\check{\text{row}}(i)\varepsilon_i,
\end{equation*}
which is a shifted half algebraic sum of positive roots for the Borel subalgebra $\mathfrak{b} \cap \hhh$ of $\hhh$.
Furthermore, we define
\begin{equation*}\label{beta}
\beta:=\sum\limits_{i\in I}(-1)^{\tp(i)}((\check{q}_1+\ldots+\check{q}_{\text{col}(i)-1})-(\check{q}_{\text{col}(i)+1}+\ldots+\check{q}_\ell))\varepsilon_i= \sum_{\alpha \in \Phi(-)}(-1)^{\tp(\alpha)} \alpha.
\end{equation*}
We introduce$$
\widetilde{\rho}:=\overline{\rho}+\beta,$$
then we have
\begin{equation}\label{tilde{rho}}
\widetilde{\rho}=\sum\limits_{i\in I}(-1)^{\tp(i)}(h-\check{\text{row}}(i)-(\check{q}_{\text{col}(i)}+\cdots+\check{q}_\ell))\varepsilon_i \
=\eta+\rho_\hhh=\frac{1}{2} \left(\sum_{\alpha \in \Phi(-) \cup \Phi(0)^+}(-1)^{\tp(\alpha)}\alpha\right) - \delta,
\end{equation}
which is a shifted half algebraic sum for the system of positive roots $\Phi(-) \cup \Phi(0)^+$.
\subsection{One-dimensional modules for $U(\ppp)$}\label{Modules for U(p)}

In this subsection, we will use the tableaux and weights defined in \textsection \ref{Tableaux} to obtain one-dimensional modules for $U(\ppp)$. Recall that we have already got the $\mathbb{Z}$-grading of $\ggg$ as in \eqref{good grading}:
$$
\ggg=\bigoplus_{k\in\mathbb{Z}}\ggg(k),\quad\text{where }  \ggg(k)=\text{span}\{e_{i,j}\in\ggg\mid \text{col}_x(j)-\text{col}_x(i)=k\}.$$
Then by \eqref{pppmmm} we have the following subalgebras of $\ggg$:
$$\ppp:=\bigoplus_{j\geq 0}\ggg(j), \qquad\qquad\hhh=\ggg(0)=\sum_{\text{col}(i)=\text{col}(j)}\Bbbk e_{i,j}.$$

By the consideration on the weights, we first give a description of the one-dimensional $U(\hhh)$-modules.
For any $i,j\in I$,
by definition  we know that 
$e_{i,j}$ is even if and only if $\tp(i)=\tp(j)$, and  $e_{i,j}$ is odd if and only if $\tp(i)\neq\tp(j)$.
Let $V:=\Bbbk v$ be a one-dimensional $U(\hhh)$-module defined by
\begin{equation}\label{eijkij}
e_{i,j}.v=k_{i,j}v
\end{equation}
with $k_{i,j}\in\Bbbk$.
Since $V$ is one-dimensional, we must have
\begin{equation}\label{kij=0}
k_{i,j}=0\qquad\text{if}~~ \tp(i)\neq\tp(j).
\end{equation}
For any $i\neq j$, we have
\begin{equation*}\label{evenlambda}
0=[e_{i,i},e_{i,j}].v=e_{i,j}.v=k_{i,j}v,
\end{equation*}
which entails that
\begin{equation}\label{kij0}
k_{i,j}=0\qquad\text{if}~~ i\neq j.
\end{equation}
It is notable that \eqref{kij0} includes \eqref{kij=0} as a special case. In fact, we further have

\begin{lemma}\label{1.2}
For a given $\lambda\in \ttt^*$, it is the weight of a one-dimensional $U(\hhh)$-module 
if and only if
\begin{equation}\label{lambdaeii}
\begin{array}{ll}
\lambda(e_{i,i})=\lambda(e_{j,j}) &\qquad\text{ whenever }\text{col}(i)=\text{col}(j) \text{ and } \tp(i)=\tp(j);\\
\lambda(e_{i,i})+\lambda(e_{j,j})=0 &\qquad\text{ whenever }col(i)=col(j) \text{ and } \tp(i)\neq \tp(j).
\end{array}
\end{equation}
\end{lemma}
\begin{proof}By the definition  we know that $e_{i,j}\in\hhh$ if and only if $\text{col}(i)=\text{col}(j)$.
Set $\lambda\in \ttt^*$ to be the weight of a one-dimensional $U(\hhh)$-module.
The discussion will proceed in two cases separately.

(1) For any $i\neq j$ with $\tp(i)=\tp(j)$, we have
\begin{equation*}\label{evenlambda}
0=[e_{i,j},e_{j,i}].v=(e_{i,i}-(-1)^{(\tp(i)+\tp(j))(\tp(i)+\tp(j))}e_{j,j}).v=(\lambda(e_{i,i})-\lambda(e_{j,j}))v,
\end{equation*}
which entails that $\lambda(e_{i,i})=\lambda(e_{j,j})$.

(2) For any $i\neq j$ with $\tp(i)\neq\tp(j)$, taking \eqref{kij=0} into consideration we have
\begin{equation*}\label{oddlambda}
0=[e_{i,j},e_{j,i}].v=(e_{i,i}-(-1)^{(\tp(i)+\tp(j))(\tp(i)+\tp(j))}e_{j,j}).v=(e_{i,i}+e_{j,j}).v=(\lambda(e_{i,i})+\lambda(e_{j,j}))v,
\end{equation*}
which entails that $\lambda(e_{i,i})+\lambda(e_{j,j})=0$.

On the other hand, for any $\lambda\in \ttt^*$ satisfies \eqref{lambdaeii}, one can easily check that $V=\Bbbk v$ defined as in \eqref{eijkij} with the restrictions \eqref{kij0} naturally becomes a $U(\hhh)$-module. We completed the proof.
%
\end{proof}

\begin{lemma}\label{1.3}Given a pyramid $\pi$,
if the $i$-th box in $\pi$ is directly above the $j$-th box and $\tp(i)=\tp(j)$, then $\widetilde{\rho}(e_{i,i})=\widetilde{\rho}(e_{j,j})+1$;
if the $i$-th box in $\pi$ is directly above the $j$-th box and $\tp(i) \neq \tp(j)$, then $\widetilde{\rho}(e_{i,i})+\widetilde{\rho}(e_{j,j})=-1.$
\end{lemma}

\begin{proof}
By virtue of $\eqref{tilde{rho}}$, we have$$
\widetilde{\rho}(e_{i,i})=(-1)^{\tp(i)}(h-\check{\text{row}}(i)-(\check{q}_{\text{col}(i)}+\cdots+\check{q}_\ell)).
$$
Assume that the $i$-th box in $\pi$ is directly above the $j$-th box. By definition we have $\text{col}(i)=\text{col}(j)$.
According to the parities of $i$ and $j$, for each case we will consider separately.
\begin{enumerate}
\item If $\tp(i)=\tp(j)=0$, then $\check{\text{row}}(j)=\check{\text{row}}(i)+1,$ so $\widetilde{\rho}(e_{i,i})=h-\check{\text{row}}(i)-(\check{q}_{\text{col}(i)}+\cdots+\check{q}_\ell)=h-\check{\text{row}}(j)+1-(\check{q}_{\text{col}(j)}+\cdots+\check{q}_\ell)=\widetilde{\rho}(e_{j,j})+1;$
\item If $\tp(i)=\tp(j)=1$, then $\check{\text{row}}(j)=\check{\text{row}}(i)-1,$ so $\widetilde{\rho}(e_{i,i})=-(h-\check{\text{row}}(i)-(\check{q}_{\text{col}(i)}+\cdots+\check{q}_\ell))=-(h-\check{\text{row}}(j)-1-(\check{q}_{\text{col}(j)}+\cdots+\check{q}_\ell))=\widetilde{\rho}(e_{j,j})+1;$
\item If $\tp(i)=0, \tp(j)=1$, then $\check{\text{row}}(j)=\check{\text{row}}(i)-1,$ so $\widetilde{\rho}(e_{i,i})=h-\check{\text{row}}(i)-(\check{q}_{\text{col}(i)}+\cdots+\check{q}_\ell)=h-\check{\text{row}}(j)-1-(\check{q}_{\text{col}(j)}+\cdots+\check{q}_\ell)=-\widetilde{\rho}(e_{j,j})-1;$
\item If $\tp(i)=1, \tp(j)=0$, then $\check{\text{row}}(j)=\check{\text{row}}(i)+1,$ so $\widetilde{\rho}(e_{i,i})=-(h-\check{\text{row}}(i)-(\check{q}_{\text{col}(i)}+\cdots+\check{q}_\ell))=-(h-\check{\text{row}}(j)+1-(\check{q}_{\text{col}(j)}+\cdots+\check{q}_\ell))=-\widetilde{\rho}(e_{j,j})-1.$
\end{enumerate}
Taking all above into consideration, the proof is completed.
\end{proof}

Now we will introduce the main result of this section.
\begin{proposition}\label{one dim to column-connected}
For $A\in \text{Tab}_\mathbb{\Bbbk}(\pi)$, $\lambda_A-\widetilde{\rho}$ is the weight of a one-dimensional $U(\hhh)$-module if and only if $A$ is column-connected.
\end{proposition}
\begin{proof}
Recall the conditions of column-connected tableau $A\in \text{Tab}_\mathbb{\Bbbk}(\pi)$  in Definition \ref{cc}.
Moreover, by the definition of $\lambda_A$ in \eqref{lambda_A} we have $\lambda_A(e_{i,i})=a_i$ for $i\in I$. By virtue of Lemma \ref{1.3}, if $A$ is column-connected we have
\begin{enumerate}
\item when the $i$-th box in $\pi$ is directly above the $j$-th box and $\tp(i)=\tp(j)$, then $\text{col}(i)=\text{col}(j)$ and $\lambda_A(e_{i,i})=\lambda_A(e_{j,j})+1$, i.e.
\begin{equation}\label{sameparity}
(\lambda_A-\widetilde{\rho})(e_{i,i})=(\lambda_A-\widetilde{\rho})(e_{j,j});
\end{equation}
\item when the $i$-th box in $\pi$ is directly above the $j$-th box and $\tp(i)\neq \tp(j)$, then $\text{col}(i)=\text{col}(j)$ and $\lambda_A(e_{i,i})+\lambda_A(e_{j,j})=-1$, i.e.
    \begin{equation}\label{diffparities}(\lambda_A-\widetilde{\rho})(e_{i,i})+(\lambda_A-\widetilde{\rho})(e_{j,j})=0.
    \end{equation}
\end{enumerate}
On the contrary, one can conclude from  \eqref{sameparity}, \eqref{diffparities} and Lemma \ref{1.3} that $A$ is column-connected.
So the proposition comes as an immediate consequence of Lemma \ref{1.2}.
%
%
%
%
%
\end{proof}

Let $\ppp':=\bigoplus_{j>0}\ggg(j)$ be a subalgebra of $\ggg$, then we have $\ppp=\hhh\oplus\ppp'$.
For a column-connected tableau $A$, we denote the above one-dimensional $U(\hhh)$-module by $\widetilde{\Bbbk}_A$.
It can be inflated to a one-dimensional $U(\ppp)$-module by the trivial $\ppp'$-action, for which we still denote by $\widetilde{\Bbbk}_A.$
\begin{remark}
In fact, for a given $A\in \text{Tab}_\mathbb{\Bbbk}(\pi)$, by the similar discussion as in Proposition \ref{one dim to column-connected} we can show that $\lambda_A - \bar \rho$ is the weight of a
one-dimensional $U(\hhh)$-module if and only if $A$ is column-connected.  If we
denote this one-dimensional $U(\hhh)$-module by $\overline \Bbbk_A$, it can also be inflated to a $U(\ppp)$-module and then induced to
a $U(\ggg)$-module $N(A) := U(\ggg) \otimes_{U(\ppp)} \overline \Bbbk_A$. It is an interesting topic to discuss when
$N(A)$ is irreducible. For more details, one can refer to \cite[Theorem 2.2]{GT} for the non-super version over an algebraically closed field of characteristic $p>0$.
\end{remark}

\section{Truncated shifted super Yangians and the commutative quotients of finite $W$-superalgebras}\label{Truncated shifted super Yangian}
In this section, we first introduce the definition of truncated shifted super Yangians, and then determine the commutative quotients of truncated shifted super Yangians and also of the corresponding finite $W$-superalgebras. All these prepare us for the final arguments on the one-dimensional representations of $U(\ggg,e)$ in \S\ref{One-dimensional modules}.

%
%
\subsection{Shifted super Yangian: Drinfeld's presentation}\label{5.1}
In this subsection we will recall the definition of shifted super Yangian $Y_{m|n}(\sigma)$, for which we refer to \cite[\S4]{P} for more details.

Let $m,n\in\mathbb{Z}_{\geq 0}$.
Recall from \S\ref{basic notations} that a pyramid $\pi$ can be uniquely determined by a triple $(\sigma,\ell,\Upsilon)$, where $\sigma$ is a shift matrix of order $m+n$ with nonnegative integral entries, $\ell$ is a positive integer and $\Upsilon$ is a fixed $0^m1^n$-sequence.
As in \cite[\S4]{P},
we use $\sigma$ and $\Upsilon$ to define the following structure, which is one of the main objects studied
in \cite{P}.


\begin{definition}\label{drshift}
Keep the notations as above.
The shifted super Yangian of $\gl_{M|N}$ associated to $\sigma$, denoted by $Y_{m|n}(\sigma)$, is the superalgebra over $\Bbbk$ generated by following symbols
\begin{align*}
&\big\lbrace D_{i}^{(r)}, D_{i}^{\prime(r)} \,|\, {1\leq i\leq m+n,\;  r\geq 0}\big\rbrace,\\
&\big\lbrace E_{j}^{(r)} \,|\, {1\leq j< m+n,\; r> s_{j,j+1}}\big\rbrace,\\
&\big\lbrace F_{j}^{(r)} \,|\, {1\leq j< m+n,\; r> s_{j+1,j}}\big\rbrace,
\end{align*}
subject to the following relations:
\begin{eqnarray}
\label{d401} D_{i}^{(0)}=D_{i}^{\prime(0)}&=&1\,,\\
\label{d402} \sum_{t=0}^{r}D_{i}^{(t)}D_{i}^{\prime (r-t)}&=&\delta_{r0},\\
\label{d403} \big[D_{i}^{(r)},D_{j}^{(s)}\big]&=&0,
\end{eqnarray}
\begin{eqnarray}
\label{d404}  [D_{i}^{(r)}, E_{j}^{(s)}]
        &=&(-1)^{\pa{i}}\big( \delta_{i,j}-\delta_{i,j+1} \big) \sum_{t=0}^{r-1} D_{i}^{(t)} E_{j}^{(r+s-1-t)}, \\
\label{d405}  [D_{i}^{(r)}, F_{j}^{(s)}]
        &=&(-1)^{\pa{i}}\big( \delta_{i,j+1} -\delta_{i,j} \big) \sum_{t=0}^{r-1} F_{j}^{(r+s-1-t)}D_{i}^{(t)}, \\
\label{d406}  [E_{i}^{(r)} , F_{j}^{(s)}]
          &=&\delta_{i,j}(-1)^{\pa{i+1}+1}
          \sum_{t=0}^{r+s-1} D_{i}^{\prime (r+s-1-t)} D_{i+1}^{(t)},
\end{eqnarray}
\begin{multline}\label{d407}
 [E_{i}^{(r)} , E_{i}^{(s)}]
          =(-1)^{\pa{i+1}}
          \big( \sum_{t=s_{i,i+1}+1}^{s-1} E_{i}^{(r+s-1-t)} E_{i}^{(t)}
          -\sum_{t=s_{i,i+1}+1}^{r-1} E_{i}^{(r+s-1-t)} E_{i}^{(t)}  \big),
\end{multline}
\begin{multline}\label{d408}
 [F_{i}^{(r)} , F_{i}^{(s)}]
          =(-1)^{\pa{i}}
          \big( \sum_{t=s_{i+1,a}+1}^{r-1} F_{i}^{(r+s-1-t)} F_{i}^{(t)}
          -\sum_{t=s_{i+1,i}+1}^{s-1} F_{i}^{(r+s-1-t)} F_{i}^{(t)}  \big),
 \end{multline}
\allowdisplaybreaks
\begin{eqnarray}
\label{d409}[E_{i}^{(r+1)}, E_{i+1}^{(s)}]-[E_{i}^{(r)}, E_{i+1}^{(s+1)}]
&=&(-1)^{\pa{i+1}} E_{i}^{(r)}E_{i+1}^{(s)}\,,\\[3mm]
\label{d410}[F_{i}^{(r+1)}, F_{i+1}^{(s)}]-[F_{i}^{(r)}, F_{i+1}^{(s+1)}]&=&
(-1)^{1+\pa{i}\pa{i+1}+\pa{i+1}\pa{i+2}+\pa{i}\pa{i+2}} F_{i}^{(s)}F_{i}^{(r)}\,,
\end{eqnarray}
\begin{align}
\label{d411}&[E_{i}^{(r)}, E_{j}^{(s)}] = 0
\qquad\qquad\text{\;\;if\;\;} |j-i|>1,\\[3mm]
\label{d412}&[F_{i}^{(r)}, F_{j}^{(s)}] = 0
\qquad\qquad\text{\;\;if\;\;} |j-i|>1,\\[3mm]
\label{d413}&\big[E_{i}^{(r)},[E_{i}^{(s)},E_{j}^{(t)}]\big]+
\big[E_{i}^{(s)},[E_{i}^{(r)},E_{j}^{(t)}]\big]=0 \quad \text{if}\,\,\, |i-j|=1,\\[3mm]
\label{d414}&\big[F_{i}^{(r)},[F_{i}^{(s)},F_{j}^{(t)}]\big]+
\big[F_{i}^{(s)},[F_{i}^{(r)},F_{j}^{(t)}]\big]=0 \quad \text{if}\,\,\, |i-j|=1,\\[3mm]
\label{d415}&\big[\,[E_{i-1}^{(r)},E_{i}^{( s_{i,i+1}+1)}]\,,\,[E_{i}^{(s_{i,i+1}+1)},E_{i+1}^{(s)}]\,\big]=0 \;\;\text{when\;\;}  m+n\geq 4\, \text{and\;\;} \pa{i}+\pa{i+1}=1,\\[3mm]
\label{d416}&\big[\,[F_{i-1}^{(r)},F_{i}^{(s_{i+1,i}+1)}]\,,\,[F_{i}^{(s_{i+1,i}+1)},F_{i+1}^{(s)}]\,\big]=0 \;\;\text{when\;\;} m+n\geq 4\, \text{and\;\;} \pa{i}+\pa{i+1}=1,
\end{align}
for all admissible indices $i,j,r,s,t$,
where their parities are given by
\begin{equation}
\label{drpa}\pa{D_{i}^{(r)}}=\pa{D_{i}^{\prime(r)}}=0, \qquad \pa{E_{j}^{(r)}}=\pa{F_{j}^{(r)}}=\pa{j}+\pa{j+1}  \,\,\text{(mod 2)}.
\end{equation}
\end{definition}

\begin{remark}
In fact, if we set
\begin{equation*}\label{Dprime}
D_{i}(u)=\textstyle{\sum}_{r\ge
0}\,D_i^{(r)}u^{-r},\quad\quad\quad D_{i}^{\prime}(u)=\textstyle{\sum}_{r\ge
0}\,D_{i}^{\prime(r)}u^{-r},
\end{equation*}
then the equation \eqref{d402} is equivalent to the condition that $D_{i}^{\prime}(u)=D_{i}(u)^{-1}$, i.e.
\begin{equation}\label{d'dre}
\textstyle{\sum}_{r\ge
0}\,D_{i}^{\prime(r)}u^{-r}\,=\,(\sum_{r\ge
0}\,D_i^{(r)}u^{-r})^{-1}=1+\sum_{k\geq1}(-1)^k(\sum_{r\ge
1}\,D_i^{(r)}u^{-r})^{k}.
\end{equation}
Then it follows from \eqref{d'dre} that $D_i^{\prime(r)}+D_i^{(r)}$ is a
polynomial in $D_i^{(1)},\ldots, D_i^{(r-1)}$ with initial form of
degree $\ge 2$; or to say, $D_i^{\prime(r)}$ can be expressed as $-D_i^{(r)}$ plus a polynomial in $D_i^{(1)},\ldots, D_i^{(r-1)}$ with initial form of
degree $\ge 2$. For example,
\begin{equation*}\label{di1pre}
D_i^{\prime(1)}=-D_i^{(1)}\quad\quad\text{and}\quad\quad
D_i^{\prime(2)}=-D_i^{(2)}+D_i^{(1)}D_i^{(1)}.
\end{equation*}
\end{remark}

%

\subsection{Commutative quotients of the truncated shifted super Yangians}\label{3.2}
In this subsection we first recall the definition of the truncated shifted super Yangian $Y_{m|n}^\ell(\sigma)$, and then determine the structure of its commutative quotient. Compared with the non-super case discussed in \cite[\S3.8]{Pre}, one can observe that the appearance of odd elements makes the situation more complicated.

We recall that $\ell=p_{m+n}$ and $p_i=\ell-s_{i,m+n}-s_{m+n,i}$ for $1\leq i\leq m+n$.
In \cite[\S8]{P}, Peng introduced the truncated shifted super Yangian of level $\ell$, denoted by $Y_{m|n}^\ell(\sigma)$, to be the quotient of $Y_{m|n}(\sigma)$ by the two-side ideal generated by the elements $\big\lbrace D_{1}^{(r)} \,|\, r> p_1 \big\rbrace$. 

For $1\leq i<j\leq m+n$, $r> s_{i,j}$ and $t>s_{j,i}$, define the following higher root elements $E_{i,j}^{(r)}, F_{j,i}^{(t)}\in Y_{m|n}(\sigma)$ recursively by
\begin{align*}
E_{i,i+1}^{(r)}:= E_{i}^{(r)}, \qquad E_{i,j}^{(r)}:=(-1)^{\pa{j-1}}[E_{i,j-1}^{(r-s_{j-1,j})}, E_{j-1}^{(s_{j-1,j}+1)}],\\
F_{i+1,i}^{(t)}:= F_{i}^{(t)}, \qquad F_{j,i}^{(t)}:=(-1)^{\pa{j-1}}[F_{j-1}^{(s_{j,j-1}+1)}, F_{j-1,i}^{(t-s_{j,j-1})}].
\end{align*}

Keep the notations as above. A PBW basis for $Y_{m|n}^\ell(\sigma)$ is obtained by Peng as follows.
\begin{lemma}(\cite[Corollary 8.3]{P})\label{spanning set}
The monomials in the elements
\begin{align*}
\big\lbrace D_{i}^{(r)} \,|\, {1\leq i\leq m+n,\;  1\leq r\leq p_i}\big\rbrace \cr \cup \big\lbrace E_{i,j}^{(r)} \,|\, {1\leq i< j\leq m+n,\; s_{i,j}<r\leq s_{i,j}+p_i}\big\rbrace \cr
\cup \big\lbrace F_{j,i}^{(r)} \,|\, {1\leq i< j\leq m+n,\; s_{j,i}<r\leq s_{j,i}+p_i} \big\rbrace
\end{align*}
taken in any fixed order span a basis for $Y_{m|n}^\ell(\sigma).$
\end{lemma}

For the $\mathbb{Z}_2$-graded superalgebra $Y_{m|n}^\ell(\sigma)=(Y_{m|n}^\ell(\sigma))_{\bar0}\oplus(Y_{m|n}^\ell(\sigma))_{\bar1}$,
let $Y_{m|n}^\ell(\sigma)^{\text{ab}}$ denote the maximal even commutative quotient of $Y_{m|n}^\ell(\sigma)$ obtained by factoring out the ideal generated by all commutators $\{[u,v]\mid u,v\in (Y_{m|n}^\ell(\sigma))_{\bar0}\}$ and all the elements in $(Y_{m|n}^\ell(\sigma))_{\bar1}$.

From now on, we will set $p_0:=0$.
The following lemma is a super version of \cite[Theorem 3.3]{Pre}. Since the existence of odd elements makes the definition of the above commutative quotients much different from the non-super case, we will give its proof in details.

\begin{proposition}\label{ab spanning set}
The commutative quotient $Y_{m|n}^\ell(\sigma)^{\text{ab}}$ is a free polynomial algebra of rank
$\ell$ generated by the elements
\begin{align}\label{ab spanning elements}
\big\lbrace d_{i}^{(r)} \,|\, {i=1,\ldots,m+n,\;  1\leq r\leq p_i-p_{i-1}}\big\rbrace
\end{align}
where $d_{i}^{(r)}$ denotes the image of $D_{i}^{(r)}$ in $Y_{m|n}^l(\sigma)^{\text{ab}}$.
\end{proposition}
\begin{proof}
Denote by  $\sfs$ the number of pairs of adjacent digits with different parities 
in $\Upsilon$. Let $j_1,j_2,\ldots,j_\sfs$ be the index of the first digit for all such pairs from left to right. Recall that the $i$-th digit of $\Upsilon$ is denoted by $|i|$ for $1\leq i\leq m+n$,
then $\Upsilon$ can be written as $(|1|,|2|,\ldots,|j_1|,|j_1+1|,\ldots,|j_2|,|j_2+1|,\ldots,|j_\sfs|,|j_\sfs+1|,\cdots,|m+n|)$, where $|j_k|$ and $|j_k+1|$ are of different parities for all $1\leq k\leq \sfs$.

(1) For the convenience of further discussion, we need the following notations.
\begin{center}
$P:=\{1,2,\ldots,m+n\}, \quad\quad S:=\{j_1,j_2,\ldots,j_\sfs,m+n\}, \quad\quad T:=P\backslash S$.
\end{center}

In virtue of \eqref{drpa}, the elements $\bigcup_{k=1}^\sfs\{E^{(r)}_{j_k}\mid r>s_{j_k,j_k+1}\}\cup\bigcup_{k=1}^\sfs\{F^{(r)}_{j_k}\mid r>s_{j_k+1,j_k}\}$  are all the odd generators of $Y_{m|n}^\ell(\sigma)$. Then the images of these elements in $Y_{m|n}^\ell(\sigma)^{\text{ab}}$ are all zeros by definition. 
Denote by $d^{(r)}_i, d_{j}^{\prime(r)} e^{(r)}_k, f^{(r)}_l$ the images of $D^{(r)}_i, D_{j}^{\prime(r)}, E^{(r)}_k, F^{(r)}_l$ in $Y_{m|n}^\ell(\sigma)^{\text{ab}}$ with $i,j\in P$ and $k,l\in T$. In particular, applying \eqref{d401} yields
\begin{equation}\label{d4012}
d^{(0)}_i=d_{j}^{\prime(0)}=1.
\end{equation}
It follows from \eqref{d402} that
\begin{equation}\label{dd'}
\sum_{t=0}^{r}d_{i}^{(t)}d_{i}^{\prime (r-t)}=\delta_{r,0}.
\end{equation}
By \eqref{d403}, the elements $D^{(r)}_i$ and $D^{(s)}_j$ commute for all admissible indices $i,j,r,s$, then we have\begin{equation}\label{dd}
[d_{i}^{(r)},d_{j}^{(s)}]=0.
\end{equation}
Applying \eqref{d404} and \eqref{d405} with $r=1$ we see that $e^{(k)}_i=f^{(k)}_i=0$ for all $i\in T$ and $k\geq1$.
On the other hand, since $[e_j^{(r)},f_j^{(1)}]=0$, if $m+n\geq2$  then \eqref{d406} entails that
\begin{equation}\label{r-tt}
          \sum_{t=0}^{r} d_{j}^{\prime (r-t)} d_{j+1}^{(t)}=0 \qquad(1\le
j \le m+n-1,\ r> p_{i+1}-p_i).
\end{equation}



Now combining all above consideration with \eqref{d401}---\eqref{d416} and the definition of $Y_{m|n}^\ell(\sigma)$, we know that $Y_{m|n}^\ell(\sigma)^{\text{ab}}$ can be described as generated by the even elements $d_i^{(k)}, d_j^{\prime(r)}$ for all admissible indices $i,j,k,r$ with the restriction $d_{1}^{(r)}=0$ for $r> p_1$, subject to the relations \eqref{d4012}---\eqref{r-tt}.

(2) We will introduce an ordinary truncated shifted Yangian, which is parallel to the algebra $Y_{m|n}^\ell(\sigma)$. Still keep $\sigma$ to be the shifted matrix of order $m+n$ with nonnegative integral entries, and the notation $\ell=p_{m+n}$ with $p_i=\ell-s_{i,m+n}-s_{m+n,i}$ for $1\leq i\leq m+n$. Just forgetting the parties of all the elements, we can define the  truncated shifted Yangian $Y_{m+n}^\ell(\sigma)$ of level $\ell$ as in \cite[\S6]{BK06}, and then its commutative quotient $Y_{m+n}^\ell(\sigma)^{\text{ab}}$. From the detailed proof as in \cite[Theorem 3.3]{Pre}, Premet showed that the algebra $Y_{m+n}^\ell(\sigma)^{\text{ab}}$ has the same symbol of generators $d_i^{(k)}, d_j^{\prime(r)}$ for all admissible indices $i,j,k,r$ with the restriction $d_{1}^{(r)}=0$ for $r> p_1$ as the ones for the algebra $Y_{m|n}^\ell(\sigma)^{\text{ab}}$, and those generators also share the same relations as in \eqref{d4012}---\eqref{r-tt} (see \cite[(2.2)\,--\,(2.15) and \S6]{BK06}). So we have $Y_{m|n}^\ell(\sigma)^{\text{ab}}\cong Y_{m+n}^\ell(\sigma)^{\text{ab}}$ as $\Bbbk$-algebras, while by \cite[Theorem 3.3]{Pre} the latter is isomorphic to a polynomial algebra in $\ell=p_{m+n}$ variables, whose generators are exactly the same symbol as in \eqref{ab spanning elements}.
\end{proof}

\subsection{Commutative quotients of the finite $W$-superalgebras}\label{W_pi}
In this subsection we will consider the commutative quotients of the finite $W$-superalgebras.


Recall the weight $\eta\in \ttt^*$ that we introduced in \eqref{eta}, which extends to a character of $\ppp$. For any $e_{i,j}\in\ppp,$ we define $$
\widetilde{e}_{i,j}:=(-1)^{\text{col}(j)-\text{col}(i)}(e_{i,j}+\eta(e_{i,j})).$$
For any $1\leq i,j\leq m+n, 0\leq x\leq m+n \text{ and } r\geq 1,$ we let
\begin{equation}\label{tdef}
T_{i,j;x}^{(r)}
:=
\sum_{s = 1}^r
\sum_{\substack{i_1,\ldots,i_s\\j_1,\ldots,j_s}}
(-1)^{\sharp\{ t=1,\ldots,s-1| \text{row}(j_t)\leq x\}}
(-1)^{\tp(i_1)+\cdots+\tp(i_s)}
 \widetilde e_{i_1,j_1} \cdots \widetilde e_{i_s,j_s}
\end{equation}
where $\sharp\{ t=1,\ldots,s-1| \text{row}(j_t)\leq x\}$ denotes the number of $t$'s with $1\leq t\leq s-1$ such that $\text{row}(j_t)\leq x$, and the sum is taken over all $i_1,\ldots,i_s,j_1,\ldots,j_s\in I$ such that
\begin{itemize}
\item[(1)] $\text{col}(j_1)-\text{col}(i_1)+\cdots+\text{col}(j_s)-\text{col}(i_s)+s= r$;
\item[(2)] $\text{col}(i_t) \leq \text{col}(j_t)$ for each $t=1,\ldots,s$;
\item[(3)] if $\text{row}(j_t)>x$, then
$\text{col}(j_t) < \text{col}(i_{t+1})$ for each
$t=1,\ldots,s-1$;
\item[(4)]
if $\text{row}(j_t)\leq x$, then $\text{col}(j_t) \geq \text{col}(i_{t+1})$
for each
$t=1,\ldots,s-1$;
\item[(5)] $\text{row}(i_1)=i$, $\text{row}(j_s) = j$;
\item[(6)]
$\text{row}(j_t)=\text{row}(i_{t+1})$ for each $t=1,\dots,s-1$.
\end{itemize}

Now define
\begin{align}
\label{d502}&D_i^{(r)}:=T_{i,i;i-1}^{(r)}\qquad\qquad\text{   for }1\leq i\leq m+n,r\geq1,\\
\label{d503}&E_i^{(r)}:=T_{i,i+1;i}^{(r)}\qquad\qquad\text{   for }1\leq i<m+n,r>s_{i,i+1},\\
\label{d504}&F_i^{(r)}:=T_{i+1,i;i}^{(r)}\qquad\qquad\text{   for }1\leq i<m+n,r>s_{i+1,i}.
\end{align}
These elements are denoted by the same symbols as the generators of the super truncated shifted Yangian and this is justified by Peng in \cite[Corollary 9.3]{P}.

Let $\mathfrak{g}=\mathfrak{gl}_{M|N}$ and $e=e_\pi\in\ggg$ be of Jordan type $(p_1,p_2,\ldots,p_{m+n})$ as in \eqref{defineepi}.
The main results of \cite{P} come as
\begin{theorem}(\cite[Theorem 10.1]{P})\label{iso}
The map from the truncated super shifted Yangian $Y_{m|n}^\ell(\sigma)$ to the finite $W$-superalgebra $U(\ggg,e)$, determined by sending each element of $Y_{m|n}^\ell(\sigma)$ in$$
\big\lbrace D_{i}^{(r)} \,|\, {1\leq i\leq m+n,\;  r\geq 1}\big\rbrace \cup \big\lbrace E_{i}^{(r)} \,|\, {1\leq i< m+n,\; r> s_{i,i+1}}\big\rbrace$$
$$\cup\big\lbrace F_{i}^{(r)} \,|\, {1\leq i< m+n,\; r> s_{i+1,i}}\big\rbrace$$
to the element of $U(\ggg,e)$ denoted by the same symbol, defines an isomorphism
$$Y_{m|n}^\ell(\sigma) \stackrel{\sim}{\rightarrow} U(\ggg,e).$$
\end{theorem}
Now we are in a position to introduce the main result of this section. By Theorem \ref{iso}, we can define $U({\mathfrak{g}},e)^{\text{ab}}$ to be the algebra $Y_{m|n}^\ell(\sigma)^{\text{ab}}$, and call it the commutative quotient of the finite $W$-superalgebra $U({\mathfrak{g}},e)$. Furthermore, from Proposition \ref{ab spanning set} we obtain the following consequence, which will play a key role in understanding the one-dimensional modules for $U(\ggg,e)$.
\begin{theorem}\label{glmn}
The commutative quotient $U({\mathfrak{g}},e)^{\text{ab}}$ is isomorphic to a polynomial algebra in $\ell=p_{m+n}$ variables \begin{align*}
\big\lbrace d_{i}^{(r)} \,|\, {i=1,\ldots,m+n,\;  1\leq r\leq p_i-p_{i-1}}\big\rbrace,
\end{align*}where the $d_{i}^{(r)}$'s have the same meaning as in Proposition \ref{ab spanning set}.
\end{theorem}

As an immediate corollary of Theorem \ref{glmn}, we have
\begin{corollary}
Set $\mathfrak{L}=\mathfrak{sl}_{M|N}$ with $M\neq N$ for $M,N\geq0$, or $\mathfrak{L}=\mathfrak{sl}(N|N)/(z)$, where $z$ is the identity matrix of order $2N$, and $(z)$ is the subalgebra of $\mathfrak{sl}(N|N)$ generated by $z$. Let  $e_\pi\in\mathfrak{L}_{\bar0}$ be of Jordan type $(p_1,p_2,\cdots,p_{m+n})$.
Then $U({\mathfrak{L}},e_\pi)^{\text{ab}}\cong\,\mathbb{C}[X_1,\ldots,X_{\ell-1}]$ with $\ell=p_{m+n}$.
\end{corollary}
\begin{proof}
Let $\ggg=\mathfrak{gl}_{M|N}$, and $e=e_\pi$.  Then we have $e\in\ggg$ and $\ggg=\mathfrak{L}\oplus\bbc z$, where $z$ is the identity matrix of order $M+N$ (for the case with $M=N$, we have $\mathfrak{gl}_{N|N}=\mathfrak{sl}_{N|N}$).
Write $\mathfrak{L}^e$ for the centralizer of $e$ in $\mathfrak{L}$. From the PBW theorem of finite $W$-superalgebras as in \cite[Theorem 4.5]{ZS1}, the operators $\Theta_i\in U(\mathfrak{L},e)$ for $1\leq i\leq \dim \mathfrak{L}^e$ given there, can be considered as elements of a generating set of the endomorphism algebra $(\text{End}_{\ggg}(U({\ggg})\otimes_{U(\mathfrak{m})}{\bbc}_\chi))^{\text{op}}=U({\ggg},e)$ corresponding to the pair $(\ggg,e)$.
Therefore, we obtain that $U(\ggg,e)\cong U(\mathfrak{L},e)\otimes_{\bbc}\bbc[z]$, which yields that
\begin{equation}\label{gablabz}
U(\ggg,e)^{\text{ab}}\cong U(\mathfrak{L},e)^{\text{ab}}\otimes_{\bbc}\bbc[z].
\end{equation}
Then $U(\mathfrak{L},e)^{\text{ab}}$ can be regarded as a subalgebra of $U(\ggg,e)^{\text{ab}}$. Taking \eqref{tdef} and \eqref{d502} into consideration, one can conclude that $z\in U(\ggg,e)$ is a linear combination of $1, D_{1}^{(1)}, \cdots, D_{m+n}^{(1)}$.

Write $\bar z$ for the image of $z$ in $U(\ggg,e)$. The preceding remark along with 
\eqref{r-tt} implies that $\bar z$ is a linear combination of $1$ and the $d^{(1)}_j$'s with $1\leq j\leq m+n$ and $p_j>p_{j-1}$. On the other hand, the proof of Proposition \ref{ab spanning set}  and then Theorem \ref{glmn} show that the algebra $U(\ggg,e)^{\text{ab}}$ is generated by $\ell$ algebraically independent elements $d^{(k)}_j$, where $1\leq j \leq m+n$, $p_j>p_{j-1}$ and $1\leq k\leq p_j-p_{j-1}$.

By all the discussions above  we know that $U(\ggg,e)^{\text{ab}}$ is a polynomial algebra in $\ell$ variables, and $\bar z$ is a linear combination of $1$ and some of the variables. It follows from \eqref{gablabz} that $\bar z\neq0$ and $U(\mathfrak{L},e)^{\text{ab}}\cong U(\ggg,e)^{\text{ab}}/(z)$. Now we can obtain that  $U({\mathfrak{L}},e)^{\text{ab}}\cong\,\mathbb{C}[X_1,\cdots,X_{\ell-1}]$, completing the proof.
\end{proof}

\section{One-dimensional modules for $U(\ggg,e)$}\label{One-dimensional modules}

%

By the similar discussion as those in \cite[\S5]{GT}, we will classify all the one-dimensional modules for $U(\ggg,e)$ in this section. As an immediate corollary of our main result stated in Theorem \ref{mean}, we will give another proof for the existence of the minimal dimensional representations of the reduced enveloping algebra $U_\chi(\ggg_\mathbb{F})$ as described in \cite[Theorem 4.3]{WZ}, where $\ggg_\mathbb{F}$ denotes the modular counterpart of $\ggg_\Bbbk$ in characteristic $p\gg0$.


\subsection{One-dimensional modules for $U(\ggg,e)$}\label{one-dimensional modules for U(g,e)}
In this subsection we will focus on the one-dimensional modules for $U(\ggg,e)$.

Let $e_r(x_1,\ldots,x_{m+n})$ denote the $r$-th elementary symmetric polynomial in the indeterminates $x_1,\ldots,x_{m+n}$. Inspired by its non-super counterpart in \cite[Lemma 2.6]{Bru}, we can obtain the following result.
\begin{lemma}\label{num}
Given a $0^m1^n$-sequence $\Upsilon$ and $a_{i}^{(r)}\in \Bbbk$ for $1\leq i\leq m+n$ and $1\leq r\leq p_i-p_{i-1},$ there exist $b_{i,j}$ for $1\leq i\leq m+n$ and $1\leq j\leq p_i$ such that
\begin{align}
\label{d604}b_{i,p_i-p_{i-1}+r}=b_{i-1,r} \qquad &\text{for} \quad 1\leq r\leq p_{i-1},\\
\label{d605}e_r((-1)^{|i|}b_{i,1},\ldots,(-1)^{|i|}b_{i,p_i})=(-1)^{r|i|}a_{i}^{(r)}\qquad &\text{for} \quad 1\leq r\leq p_i-p_{i-1}.
\end{align}
\end{lemma}
\begin{proof}
We prove the existence of the numbers $b_{i,j}$ for $1\leq j\leq p_i$ satisfying the above relation by induction on $i=1,\ldots,m+n$.

(1) Consider the basic case $i=1$.

(i) When $|1|=0$, the equation \eqref{d605} becomes $e_r(b_{1,1},\ldots,b_{1,p_1})=a_{1}^{(r)}$. Thus we define $b_{1,1},\ldots,b_{1,p_1}\in \Bbbk$ from the factorization
$$
t^{p_1}+a_1^{(1)}t^{p_1-1}+\cdots+a_1^{(p_1)}=(t+b_{1,1})\cdots(t+b_{1,p_1}).
$$

(ii) When $|1|=1,$ the equation \eqref{d605} is $e_r(-b_{1,1},\cdots,-b_{1,p_1})=(-1)^{r}a_{1}^{(r)}$. Thus we  define $b_{1,1},\ldots,b_{1,p_1}\in \Bbbk$ from the factorization
$$
t^{p_1}-a_1^{(1)}t^{p_1-1}+\cdots+(-1)^{p_1}a_1^{(p_1)}=(t-b_{1,1})\cdots(t-b_{1,p_1}).
$$

(2) For the induction step, suppose that we have already defined $b_{i-1,1},\ldots,b_{i-1,p_{i-1}}$. First define $b_{i,p_i-p_{i-1}+1},\ldots,b_{i,p_i}$ so that the equation \eqref{d604} holds. Then we need to find numbers $b_{i,1},\ldots,b_{i,p_i-p_{i-1}}$ satisfying \eqref{d605}.

(i) If $|i|=0,$ the equations \eqref{d605} are equivalent to the equations
\begin{align*}
e_r(b_{i,1},\ldots,b_{i,p_i-p_{i-1}})&=a_{i}^{(r)}-\sum_{s=0}^{r-1}e_s(b_{i,1},\ldots,b_{i,p_i-p_{i-1}})e_{r-s}(b_{i,p_i-p_{i-1}+1},\ldots,b_{i,p_i})\\
&=a_{i}^{(r)}-\sum_{s=0}^{r-1}e_s(b_{i,1},\ldots,b_{i,p_i-p_{i-1}})e_{r-s}(b_{i-1,1},\ldots,b_{i-1,p_{i-1}}).
\end{align*}
for $1\leq r\leq p_i-p_{i-1}.$ Carrying induction on $r=1,\ldots,p_i-p_{i-1},$ we solve the equations uniquely for $b_i^{(r)}:=e_r(b_{i,1},\ldots,b_{i,p_i-p_{i-1}})$, and then define $b_{i,1},\ldots,,b_{i,p_i-p_{i-1}}$ by factoring
$$
t^{p_i-p_{i-1}}+b_i^{(1)}t^{p_i-p_{i-1}-1}+\cdots+b_i^{(p_i-p_{i-1})}=(t+b_{i,1})\cdots(t+b_{i,p_i-p_{i-1}}).
$$

(ii) If $|i|=1,$ the equations \eqref{d605} are equivalent to the equations
\begin{align*}
e_r(-b_{i,1},\ldots,-b_{i,p_i-p_{i-1}})&=(-1)^{r}a_{i}^{(r)}-\sum_{s=0}^{r-1}e_s(-b_{i,1},\ldots,-b_{i,p_i-p_{i-1}})e_{r-s}(-b_{i,p_i-p_{i-1}+1},\ldots,-b_{i,p_i})\\
&=(-1)^{r}a_{i}^{(r)}-\sum_{s=0}^{r-1}e_s(-b_{i,1},\ldots,-b_{i,p_i-p_{i-1}})e_{r-s}(-b_{i-1,1},\ldots,-b_{i-1,p_{i-1}})
\end{align*}
for $1\leq r\leq p_i-p_{i-1}.$ Proceeding by induction on $r=1,\ldots,p_i-p_{i-1},$ we can obtain the numbers $b_i^{\prime(r)}:=e_r(-b_{i,1},\ldots,-b_{i,p_i-p_{i-1}})$, and then define $b_{i,1},\ldots,b_{i,p_i-p_{i-1}}$ by factoring
$$
t^{p_i-p_{i-1}}+b_i^{\prime(1)}t^{p_i-p_{i-1}-1}+\cdots+b_i^{\prime(p_i-p_{i-1})}=(t-b_{i,1})\cdots(t-b_{i,p_i-p_{i-1}}).
$$Then the lemma follows.
\end{proof}

Due to the property of the elementary symmetric polynomials, we can permute the 
numbers $b_{i,j}$ as in Lemma \ref{num} for $1\leq j\leq p_i$   with each fixed $i$.
Then we have the following lemma.
\begin{lemma}\label{cn}
Given a $0^m1^n$-sequence $\Upsilon$, a shift matrix $\sigma=(s_{i,j})_{(m+n)\times(m+n)}$ and $a_{i}^{(r)}\in \Bbbk$ for $i=1,\ldots,m+n$ and $1\leq r\leq p_i-p_{i-1},$ there exist $b_{i,j}$ for $1\leq i\leq m+n$ and $1\leq j\leq p_i$ such that
\begin{align}
\label{d604'}b_{i,s_{i,i-1}+r}=b_{i-1,r} \qquad &\text{for} \quad 1\leq r\leq p_{i-1},\\
\label{d605'}e_r((-1)^{|i|}b_{i,1},\ldots,(-1)^{|i|}b_{i,p_i})=(-1)^{r|i|}a_{i}^{(r)}\qquad &\text{for} \quad 1\leq r\leq p_i-p_{i-1}.
\end{align}
\end{lemma}

Now  we are in a position to give the classification of one-dimensional $U(\ggg,e)$-modules. To begin with, we first make some conventions. 
For any pyramid $\pi$, we know that all the boxes in the same row have the same ``$+$" or ``$-$" labeling. Let $j_1<\ldots<j_{p_i}$ in $I$ be all the elements in the $i$-th row counting from left to right. Then we have $|i|=\tp(j_1)=\ldots=\tp(j_{p_i})$, and also $\check{\text{row}}(j_1)=\ldots=\check{\text{row}}(j_{p_i})$ by definition. 
To simply notations, we will write $\hat{\text{row}}(i):=\check{\text{row}}(j_1)=\ldots=\check{\text{row}}(j_{p_i})$ in the following discussion.

For $i=1,\ldots,m+n,$ we write $a_{i,1},\ldots,a_{i,p_i}$ for the entries in the $i$-th row of $A\in \text{Tab}_\Bbbk(\pi)$ from left to right. It follows from (\ref{d403}) and Lemma \ref{spanning set}   that $\big\lbrace D_{i}^{(r)} \,|\, {1\leq i\leq m+n,\;  1\leq r\leq p_i}\big\rbrace$ generates a subalgebra of $U(\ggg,e)$ isomorphic to a polynomial algebra in $m+n$ variables. We denote this subalgebra by $U(\ggg,e)^0$. We define the one-dimensional $U(\ggg,e)^0$-module $\Bbbk_{\overline{A}}$ by letting $D_i^{(r)}$ act on $\Bbbk_{\overline{A}}$ by $(-1)^{r|i|}e_r(a_{i,1}+(-1)^{|i|}\hat{\text{row}}(i),\ldots,a_{i,p_i}+(-1)^{|i|}\hat{\text{row}}(i))$. 
It is obvious that $\Bbbk_{\overline{A}}$ depends only the row equivalence class of $A$. We note that, for any given $b_1,\ldots,b_{p_i}\in \Bbbk,$ finding $c_{i,1},\ldots,c_{i,p_i}$ such that $e_r(c_{i,1},\ldots,c_{i,p_i})=b_r$ for each $r$ is equivalent to find solutions of the polynomial $t^{p_i}-b_1t^{p_i-1}+\cdots+(-1)^{p_i}b_{p_i}.$ Therefore, since $\Bbbk$ is algebraically closed, we see that any one-dimensional $U(\ggg,e)^0$-module is isomorphic to $\Bbbk_{\overline{A}}$ for some $A\in \text{Tab}_\Bbbk(\pi).$ Thus we see that the restriction of any one-dimensional $U(\ggg,e)$-module to $U(\ggg,e)^0$ is isomorphic to $\Bbbk_{\overline{A}}$ for some $A\in \text{Tab}_\Bbbk(\pi),$ and that such $A$ is defined up to row equivalence. By the same consideration as Step (1) in the proof of Proposition \ref{ab spanning set}, we see that the generators $E_i^{(s)}$ and $F_i^{(s)}$ act as $0$ on any one-dimensional $U(\ggg,e)$-module, for all $i$ and $s$. If the action of $U(\ggg,e)^0$ on $\Bbbk_{\overline{A}}$ can be extended to a $U(\ggg,e)$-module, on which $E_i^{(s)}$ and $F_i^{(s)}$ act as 0, then we denote this module by $\widetilde{\Bbbk}_{\overline{A}}$. We will determine when $\widetilde{\Bbbk}_{\overline{A}}$ exists, and this is achieved in the following theorem.
\begin{theorem}\label{mean}
For any $A\in\text{Tab}_\Bbbk(\pi)$, there is a one-dimensional $U(\ggg,e)$-module $\widetilde{\Bbbk}_{\overline{A}}$, which extends the action of $U(\ggg,e)^0$ on $\Bbbk_{\overline{A}}$ if and only if $A$ is row-equivalent to a column-connected tableau.
\end{theorem}
\begin{proof}
(1) First let $A\in \text{Tab}_\Bbbk(\pi)$ be column-connected with entries in the $i$-th row labelled $a_{i,1},\ldots,a_{i,p_i}$ for $i=1,\ldots,m+n.$ Recall the one-dimensional $U(\ppp)$-module $\widetilde{\Bbbk}_A$ we introduced in \S\ref{Modules for U(p)}. Consider the action of the explicit elements $D_{i}^{(r)}\in U(\ppp)$ given in (\ref{tdef}) and (\ref{d502}) on the module $\widetilde{\Bbbk}_A$. The only summands in the expression for $D_{i}^{(r)}$ that do not act as zero on $\widetilde{\Bbbk}_A$ are those which are $(-1)^{r|i|} \widetilde e_{j_1,k_1} \cdots \widetilde e_{j_s,k_s}$ such that $j_1=k_1,\ldots,j_s=k_s$ and $\text{row}(j_1)=\ldots=\text{row}(j_s)=i$: terms of this form only occur for $s=r$ and their sum is $(-1)^{r|i|}\widetilde e_{j_1,j_1} \cdots \widetilde e_{j_r,j_r}$ for all $j_1<j_2<\cdots<j_r$ of $I$  in the $i$-th row. The $\mathfrak{t}$-weight of $\widetilde{\Bbbk}_A$ is $\lambda_A-\widetilde{\rho}$, and we have $\widetilde{\rho}=\eta+\rho_\mathfrak{h}.$ Thus we see that each $\widetilde e_{j,j}$ acts on $\widetilde{\Bbbk}_A$ by $(\lambda_A-\rho_\mathfrak{h})(e_{j,j}),$ because of the shift of $\eta$ in the definition of $\widetilde e_{j,j}$. From all of these observations we know that the action of $D_{i}^{(r)}$ on $\widetilde{\Bbbk}_A$ is given by $(-1)^{r|i|}e_r(a_{i,1}+(-1)^{|i|}\hat{\text{row}}(i),\ldots,a_{i,p_i}+(-1)^{|i|}\hat{\text{row}}(i))$. 
This proves the existence of $\widetilde{\Bbbk}_{\overline{A}}$ under the assumption that $A$ is column-connected.

(2)
We move on to prove that $\widetilde{\Bbbk}_{\overline{A}}$ exists only if $A$ is column-connected. First note that the action $U(\ggg,e)$ on any one-dimensional module factors to an action of the commutative quotient $U(\ggg,e)^{\text{ab}}$ of $U(\ggg,e).$ Thanks to Theorem \ref{glmn}, we know that the algebra $U(\ggg,e)^{\text{ab}}$ is generated by the images of the $\ell$ elements$$
\big\lbrace D_{i}^{(r)} \,|\, {i=1,\ldots,m+n,\;  0< r\leq p_i-p_{i-1}}\big\rbrace.$$
So any one-dimensional $U(\ggg,e)$-module is determined uniquely by the action of these elements. Thus to show that any one-dimensional $U(\ggg,e)$-module is of the form $\widetilde{\Bbbk}_{\overline{A}}$ for some column-connected $A\in  \text{Tab}_\Bbbk(\pi)$, it is suffices to show that for any set $$
\big\lbrace a_{i}^{(r)} \,|\, {i=1,\ldots,m+n,\;  0< r\leq p_i-p_{i-1}}\big\rbrace$$
with $a_{i}^{(r)}\in \Bbbk,$ there is a column-connected $A\in  \text{Tab}_\Bbbk(\pi)$ such that the action of $D_{i}^{(r)}$ on $\widetilde{\Bbbk}_{\overline{A}}$ is given by $a_{i}^{(r)}$ for $i=1,\ldots,m+n$ and $0< r\leq p_i-p_{i-1}.$

In fact,
by lemma \ref{cn}  we can find $b_{i,j}$ for $1\leq i\leq m+n$ and $1\leq j\leq p_i$ such that
\begin{align*}
b_{i,s_{i,i-1}+r}=b_{i-1,r} \qquad &\text{for} \quad 1\leq r\leq p_{i-1},\\
e_r((-1)^{|i|}b_{i,1},\ldots,(-1)^{|i|}b_{i,p_i})=(-1)^{r|i|}a_{i}^{(r)}\qquad &\text{for} \quad 1\leq r\leq p_i-p_{i-1}.
\end{align*}
So we can choose $$
a_{i,j}:=(-1)^{|i|}(b_{i,j}-\hat{\text{row}}(i))$$ for $1\leq j \leq p_i$,
and take $a_{i,1},\ldots,a_{i,p_i}$ as the entries of the $i$-th row from left to right, which form a tableau $A\in \text{Tab}_\Bbbk(\pi).$ By the same consideration as in the proof of Lemma \ref{1.3}, we can prove that
%
$A$ is column-connected. Recall in Step (1) we have shown that $D_{i}^{(r)}$ acts on $\widetilde{\Bbbk}_A$ by
$(-1)^{r|i|}e_r(\widetilde e_{i_1,i_1} \cdots \widetilde e_{i_r,i_r})$.
Then we can define the action of $D_{i}^{(r)}$ on $\widetilde{\Bbbk}_A$ by
$(-1)^{r|i|}e_r(a_{i,1}+(-1)^{|i|}\hat{\text{row}}(i),\ldots,a_{i,p_i}+(-1)^{|i|}\hat{\text{row}}(i))=a_i^{(r)}$ for $i=1,\ldots,m+n$ and $0< r\leq p_i-p_{i-1}$.
From this we can deduce that all one-dimensional $U(\ggg,e)$-modules are of the form $\widetilde{\Bbbk}_{\overline{A}}$ for some column-connected $A\in\text{Tab}_\Bbbk(\pi)$.
\end{proof}

Thanks to Theorem \ref{mean}, we can give the complete classification of one-dimensional $U(\ggg,e)$-modules.
\begin{corollary}\label{meancoro}
As $A\in\text{Tab}_\Bbbk(\pi)$ runs over a set of representations for the row-equivalent classes of column-connected tableaux, the set of the modules $\widetilde{\Bbbk}_{\overline{A}}$ gives a complete set of pairwise-inequivalent one-dimensional $U(\ggg,e)$-modules.
\end{corollary}
\begin{remark}
Similar to the non-super case, it is worth to mention that the freeness of the generators $d_{i}^{(r)}$'s for the commutative quotient $Y_{m|n}^\ell(\sigma)^{\text{ab}}$ as in Proposition \ref{ab spanning set} can also be deduced from Step (2) in the proof of Theorem \ref{mean}, which shows that there are one-dimensional modules for $Y_{m|n}^\ell(\sigma)\cong U(\ggg,e)$ on which these generators
can act by arbitrary elements of $\Bbbk$, thus skipping the discussion on the structure of commutative quotient of the ordinary  truncated shifted Yangian $Y_{m+n}^\ell(\sigma)$ as in Step (2) for the proof of Proposition \ref{ab spanning set}.
\end{remark}

\subsection{On the existence of the minimal dimensional modules for $\gl_{M|N}$ in characteristic $p\gg0$}
In this subsection we will consider the minimal dimensional modules for $\gl_{M|N}$ in characteristic $p\gg0$.

Let $\ggg$ be a basic classical Lie superalgebra over $\Bbbk$ and $e$ be an even nilpotent element in $\ggg$. Set $(\cdot,\cdot)$ to be an even non-degenerate supersymmetric invariant bilinear form on $\ggg$, and define $\chi\in{\ggg}^{*}$ by letting $\chi(x)=(e,x)$ for all $x\in{\ggg}$.
Denote by $\ggg^e$ the centralizer of $e$ in $\ggg$. Write $d_0:=\text{dim}\,\ggg_{\bar 0}-\text{dim}\,\ggg^e_{\bar 0}$ and $d_1:=\text{dim}\,\ggg_{\bar 1}-\text{dim}\,\ggg^e_{\bar 1}$.

Let $\ggg_{\mathbb{F}}$ be the modular counterpart of $\ggg$, and still write $\chi\in({\ggg}_{\mathbb{F}})^*_{\bar0}$ for the modular version of $\chi\in{\ggg}^{*}$ as defined above. Then $\chi$ is a nilpotent $p$-character. Denote by $U_\chi({\ggg}_{\mathbb{F}})$ the reduced enveloping algebra of $\ggg_{\mathbb{F}}$ associated to $\chi$. Recall that any irreducible ${\ggg}_{\mathbb{F}}$-module $V$ is a module of $U_\xi({\ggg}_{\mathbb{F}})$ for a unique $\xi=\xi_V\in({\ggg}_{\mathbb{F}})_{\bar{0}}^*$. We often regard $\xi\in{\ggg}_{\mathbb{F}}^*$ by letting $\xi(({\ggg}_{\mathbb{F}})_{\bar{1}})=0$.

In \cite[Theorem 4.3]{WZ}, Wang-Zhao introduced the super Kac–Weisfeiler property, which says that under the restriction imposed on the characteristic of the field $\mathbb{F}$ as in \cite[Table 1]{WZ}, the dimension of any irreducible representation of $U_\chi({\ggg}_{\mathbb{F}})$ is divisible by the number $p^{\frac{d_0}{2}}2^{\lfloor\frac{d_1}{2}\rfloor}$, where $\lfloor a\rfloor$ denotes the least integer upper bound of $\frac{d_1}{2}$.
With this number, Zeng-Shu partially answered the question whether there exist modules whose dimensions is equal to such a number in  \cite{ZS4}, i.e.,
\begin{theorem}(\cite[Theorem 1.6]{ZS4})\label{intromain-2}
Let ${\ggg}_\mathbb{F}$ be a  basic classical Lie superalgebra over $\mathbb{F}=\overline{\mathbb{F}}_p$, and let $\chi\in({\ggg}_{\mathbb{F}})^*_{\bar0}$ be a nilpotent $p$-character, with respect to the element $e\in(\ggg_\mathbb{F})_\bz$ under the bilinear form $(\cdot,\cdot)$ as above. If the corresponding finite $W$-superalgebra $U({\ggg},e)$ over ${\bbc}$ affords a one-dimensional  (resp. two-dimensional) representation when $d_1$ is even (resp. odd), then for $p\gg0$ the reduced enveloping algebra $U_\chi({\ggg}_{\mathbb{F}})$ admits irreducible representations of dimension $p^{\frac{d_0}{2}}2^{\lfloor\frac{d_1}{2}\rfloor}$.
\end{theorem}


In particular, for the case $\ggg=\gl_{M|N}$ we know that $d_1$ is always even by \cite[\S3.2]{WZ}. Applying the tool of the so-called ``shifted pyramids", Zeng-Shu further gave the following result. 
\begin{theorem}(\cite[Proposition 3.1]{ZS2})\label{main2}
Let $\ggg_\mathbb{F}=\mathfrak{gl}_{M|N}(\mathbb{F})$ over a field $\mathbb{F}$ of characteristic $p>2$. For any nilpotent $p$-character $\chi\in\ggg^*_{\bar0}$, the reduced enveloping algebra $U_\chi(\ggg_\mathbb{F})$ admits irreducible modules of dimension $p^{\frac{d_0}{2}}2^{\frac{d_1}{2}}$.
\end{theorem}

By the method of finite $W$-superalgebras, we can give another proof of Corollary \ref{mincor}, which is the weak version of Theorem \ref{main2} in characteristic $p\gg0$. In fact, Corollary \ref{mincor} can be easily deduce from 
Theorems \ref{mean} and \ref{intromain-2}.
\begin{remark}
Although the consequence in Corollary \ref{mincor} seems to be weaker than the origin one Zeng-Shu obtained in \cite{ZS2}, it takes the first step towards the profound understanding of $U_\chi(\mathfrak{gl}_{M|N}(\mathbb{F}))$-modules of minimal dimensions. In fact, for the case with $\mathfrak{l}_\mathbb{F}=\gl_N(\mathbb{F})$ in any characteristic $p>0$,  Goodwin-Topley obtained all the isomorphism class of $U_\chi(\mathfrak{l}_\mathbb{F})$-modules of minimal dimensions in \cite{GT}, where the complete classification of one-dimensional representations of the finite $W$-algebra $U(\mathfrak{l}_\mathbb{F},e)$ 
plays a critical role. It is hopeful to get similar results for $\mathfrak{gl}_{M|N}(\mathbb{F})$ by generalizing the consequences we obtained in the present paper to their modular counterparts.
\end{remark}

\subsection*{Acknowledgements}
Both authors would like to express sincere gratitude to their supervisor, Professor Bin Shu, for his constant encouragement and guidance. The authors also would like to thank Hao Chang and Yung-Ning Peng for stimulating discussions.

\end{document}